\documentclass{amsart}
\usepackage[pdfauthor={Richard J. Mathar},pdfkeywords={Gauss Quadrature, Analysis, Logarithmic Kernel},pdfsubject={Mathematics, Analysis},colorlinks=true,citecolor=blue]{hyperref}

\usepackage{amsmath}

\newtheorem{rem}{Remark}

\begin{document}

\title{Gaussian Quadrature of $\int_0^1 f(x) \log^m(x) dx$ and $\int_{-1}^1 f(x) \cos(\pi x/2)dx$}

\author{Richard J. Mathar} 
\email{mathar@mpia.de}
\urladdr{http://www.mpia.de/~mathar}
\address{Hoeschstr. 7, 52372 Kreuzau, Germany}


\date{\today}

\begin{abstract}
We tabulate the abscissae and associated weights for numerical integration of
integrals with either the singular weight function $(-\log x)^m$
for exponents $m=1,2$ or $3$,
or the symmetric weight function $\cos(\pi x/2)$.
Standard brute force arithmetics generates explicit pairs of these
values for up to 128 nodes.
\end{abstract}

\maketitle 

\section{Methodology} 

The paper provides abscissae $x_i$ and weights $w_i$ for
Gaussian integration with a power of a logarithm
in the integral kernel on one hand,
\begin{equation}
\int_0^1 f(x) (-\log x)^m dx
\approx \sum_{i=1}^N w_i f(x_i).
\label{eq.log}
\end{equation}
or with a cosine in the integral kernel on the other,
\begin{equation}
\int_{-1}^1 f(x) \cos(\pi x/2) dx
\approx \sum_{i=1}^N w_i f(x_i).
\label{eq.cos}
\end{equation}
The $w_i$ and $x_i$ are computed with the
standard theory from roots of a system of orthogonal polynomials $p_n$
with norm
\cite{Golub,Luke1975,WynnQJM18}
\begin{equation}
\langle f,g \rangle \equiv \int_0^1 f(x) g(x)(-\log x)^m dx,
\label{eq.logortho}
\end{equation}
and
\begin{equation}
\langle f,g \rangle \equiv \int_{-1}^1 f(x) g(x)\cos \frac{\pi x}{2} dx,
\end{equation}
respectively.
A set of orthogonal (monic) polynomials $p_{n}(x)$ is bootstrapped from
\begin{equation}
p_{-1}(x)=0;\quad  p_{0}(x)=1; \quad p_{n+1}(x)= (x-a_{n})p_{n}(x)-b_{n} p_{n-1}(x).
\label{eq.rec}
\end{equation}
[Dependence of polynomials and coefficients $a$ and $b$ on the parameter $m$
in the case (\ref{eq.log}) is not written down explicitly here.] Multiplication of
the recurrence with $p_n$ or $p_{n-1}$ and using the requirement of orthogonality
proposes to calculate the coefficients and polynomials recursively with
\begin{equation}
a_{n} = \frac{\langle xp_{n},p_{n}\rangle}{\langle p_{n},p_{n}\rangle};
\end{equation}
\begin{equation}
b_0=0;\quad b_{n} = \frac{\langle xp_{n},p_{n-1}\rangle}{\langle p_{n-1},p_{n-1}\rangle} \quad (n>0).
\label{eq.bofn}
\end{equation}
\begin{rem}
In cases like (\ref{eq.cos}) where the weight in the integral is an even function and the integral
limits are symmetric, all $a_n$ are zero.
\end{rem}
The standard further steps are
\begin{itemize}
\item
normalization of the polynomials such that
their norm is unity,
\begin{equation}
p^*_n(x) \equiv \frac{p_n(x)}{\sqrt{\langle p_n,p_n\rangle}},
\end{equation}
\item
computation of all zeros $x_i$ of $p_N(x)$ at some degree $N$.
\item
computation of the weights $w_i$ by
\begin{equation}
w_i= -\frac{[x^{N+1}]p_{N+1}^*}{ [x^N]p_N^*}\, \frac{1}{ p_{N+1}^*(x_i){p_N^*}'(x_i)},
\end{equation}
where $[x^{N+1}]p_{N+1}^*$ and $[x^N]p_N^*$ are the leading coefficients of the
two polynomials after normalization, and where the prime at $p'$ denotes the derivative
with respect to $x$.
\end{itemize}
We obviously add no new aspect to the established theory.
The benefit is to those readers who need explicit abscissae-weight pairs
and have no access to a multi-precision numeric library.

\section{Logarithmic Kernel} 

The first part of the results extends tables that have been published
in the literature for exponent $m=1$, namely by Anderson
for $N$ up to 10 \cite{AndersonMC19},
by Danloy for $N=10$ and $N=20$ \cite{DanloyMC27},
and by King for $N=20$ and $N=30$ \cite{KingJOSAB19}.

\begin{rem}
The variable substitution $x= e^{-y}$
changes the format to
\begin{equation}
\int_0^1 f(x)(-\log x)^m dx
=
\int_0^\infty f(e^{-y})y^m e^{-y} dy
\end{equation}
which is alternatively evaluated with
Gauss-Laguerre quadratures
\cite[(25.4.38)]{AS}\cite{ChisholmJCP10,SalzerCJ15}.
\end{rem}
Integrals of the form (\ref{eq.logortho}) are calculated for the polynomials that appear
in the recurrence (\ref{eq.rec}) term-by-term with the aid of the moments
$\mu$ \cite[2.722]{GR},
\begin{equation}
\mu_{n,m}\equiv \int_0^1 x^n (-\log x)^m dx
=\frac{m!}{(n+1)^{m+1}}
.
\end{equation}

The first polynomials $p_{n,m}(x)$ look as follows:
\begin{eqnarray}
p_{1,1} &=& x-1/4 ;\\
p_{2,1} &=& {x}^{2}-5/7\,x+{\frac {17}{252}} ;\\
p_{3,1} &=& {x}^{3}-{\frac {3105}{2588}}\,{x}^{2}+{\frac {5751}{16175}}\,x-{\frac {
4679}{258800}} ;\\
p_{1,2} &=& x-1/8 ;\\
p_{2,2} &=& {x}^{2}-{\frac {19}{37}}\,x+{\frac {217}{7992}} ;\\
p_{3,2} &=& {x}^{3}-{\frac {1632663}{1695176}}\,{x}^{2}+{\frac {5619807}{26487125}} \,x-{\frac {1568083}{242168000}} ;\\
p_{1,3} &=& x-1/16 ;\\
p_{2,3} &=& {x}^{2}-{\frac {13}{35}}\,x+{\frac {493}{45360}}; \\
p_{3,3} &=& {x}^{3}-{\frac {129197997}{166534960}}\,{x}^{2}+{\frac {4147011999}{
32526359375}}\,x-{\frac {19126701359}{8326748000000}}.
\end{eqnarray}
$p_{2,1}$ in particular has been written down earlier \cite{LeanIJNM21}.
Two generic values are
\begin{eqnarray}
p_{1,m} &=&
x-{2}^{-1-m}; \\
p_{2,m} &=&
x^2
+
\frac {-2^{m+1}+3^{m+1}}{3^{m+1}-4^{m+1}}x
+
\frac {-3^{m+1}{2}^{-1-m}+4^{m+1}{3}^{-1-m}}{3^{m+1}-4^{m+1}}.
\end{eqnarray}

The results are summarized in the ASCII files \texttt{log\_}$N$\_$m$ in the ancillary
directory, where $N$ covers the range $3$ to $128$ and $m$ covers powers from $1$ to $3$.
Each line contains a pair $(x_i,w_i)$.
For improved readability, a blank line is inserted after each block of 5 nodes.
The numbers have been stabilized to the 30 digits shown by cranking up
the internal representation of numbers in a Maple program to 270 digits.

\begin{rem}
Related approximative cubatures where polynomials are not only multiplied by
also added to the logarithm in the kernel have also been discussed
\cite{HarrisIJCM6,CrowMC60,KolmCMA41}.
\end{rem}

\section{Cosine kernel}

The tools to assemble (\ref{eq.cos}) start
from repeated partial integration of
\cite[3.761]{GR}
\begin{equation}
\int_0^{\pi/2} x^m\cos x dx = \sum_{k=0}^{\lfloor m/2\rfloor}
(-)^k\frac{m!}{(m-2k)!}\left(\frac{\pi}{2}\right)^{m-2k}
+(-)^{\lfloor m/2\rfloor}m!(2\lfloor \frac{m}{2}\rfloor-m)
,
\end{equation}
for non-negative integer $m$.
The even moments are therefore
\begin{multline}
\mu_{2m}\equiv \int_{-1}^1 x^{2m}\cos (x\pi/2) dx
=
2(2m)!\sum_{k=0}^m
(-)^k\frac{1}{(2m-2k)!}\left(2/\pi\right)^{2k+1}
\\
=
\frac{4}{\pi} \,_3F_0\left(\begin{array}{c} -m+\frac{1}{2},-m,1\\ - \end{array}\mid -\frac{16}{\pi^2}\right)
.
\end{multline}
The odd moments are zero because the cosine is an even function.
The monic orthogonal polynomials start
\begin{eqnarray}
p_0 &=&1;\quad p_1=x;\\
p_2 &=& x^2-1+\frac{8}{\pi^2};\\
p_3 &=& x^3-\frac{\pi^4-48\pi^2+384}{(\pi^2-8)\pi^2}x;\\
p_4 &=& x^4-2\frac{\pi^4-78\pi^2+672}{\pi^2(\pi^2-10)}x^2
+\frac{\pi^6-114\pi^4+1728\pi^2-6912}{\pi^4(\pi^2-10)},
\end{eqnarray}
and have parities $p_{-n}(x)=(-1)^np_n(x)$.

The results are summarized in the ASCII files \texttt{cosine\_}$N$ in the ancillary
directory, where $N$ covers the range $3$ to $128$.
The numbers have been stabilized to the 30 digits shown by
an internal representation of numbers in a Maple program with 650 digits.

Only the values with positive $x_i$ or $x_i=0$ are tabulated; the
duplicates of the nodes at the negative abscissae (with the same weights) are not
added explicitly.

\bibliographystyle{amsplain}
\bibliography{all}

\providecommand{\bysame}{\leavevmode\hbox to3em{\hrulefill}\thinspace}
\providecommand{\MR}{\relax\ifhmode\unskip\space\fi MR }
\providecommand{\MRhref}[2]{%
  \href{http://www.ams.org/mathscinet-getitem?mr=#1}{#2}
}
\providecommand{\href}[2]{#2}
\begin{thebibliography}{10}

\bibitem{AS}
Milton Abramowitz and Irene~A. Stegun (eds.), \emph{Handbook of mathematical
  functions}, 9th ed., Dover Publications, New York, 1972. \MR{0167642 (29
  \#4914)}

\bibitem{AndersonMC19}
Donald~G. Anderson, \emph{Gaussian quadrature for $\int_0^1 -ln(x) f(x) dx$},
  Math. Comput. \textbf{19} (1965), 477--481. \MR{0178569}

\bibitem{ChisholmJCP10}
J.~S.~R. Chisholm and A.~Genz, \emph{Accelerated convergence of sequences of
  quadrature approximations}, J. Comput. Phys. \textbf{10} (1972), no.~2,
  284--307. \MR{0326998}

\bibitem{CrowMC60}
John~A. Crow, \emph{Quadrature of integrands with a logarithmic singularity},
  Math. Comput. \textbf{60} (1993), no.~201, 297--301.

\bibitem{DanloyMC27}
Bernard Danloy, \emph{Numerical construction of gaussian quadrature formuals
  for $\int_0^1(-\log x) x^\alpha f(x)dx$ and $\int_0^\infty e_m(x) f(x) dx$},
  Math. Comp. \textbf{27} (1973), no.~124, 861--869. \MR{0331730}

\bibitem{Golub}
Gene~H. Golub and John~H. Welsch, \emph{Calculation of {G}auss {Q}uadrature
  {R}ules}, Math.\ Comp. \textbf{23} (1969), no.~106, 221--230. \MR{0245201 (39
  \#6513)}

\bibitem{GR}
I.~Gradstein and I.~Ryshik, \emph{Summen-, {P}rodukt- und {I}ntegraltafeln},
  1st ed., Harri Deutsch, Thun, 1981. \MR{0671418 (83i:00012)}

\bibitem{HarrisIJCM6}
C.~G. Harris and W.~A.~B. Evans, \emph{Extension of numerical quadrature
  formuale to cater for end point singuar behaviours over finite intervals},
  Int. J. Computer Math. B \textbf{6} (1977), 219--227.

\bibitem{KingJOSAB19}
Frederick~W. King, \emph{Efficient numerical approach to the evaluation of
  {K}ramers-{K}ronig transforms}, J. Opt. Soc. Am. B \textbf{19} (2002),
  no.~10, 2427--2436. \MR{1945703}

\bibitem{KolmCMA41}
P.~Kolm and V.~Rokhlin, \emph{Numerical quadratures for singular and
  hypersingular integrals}, Comp. Math. Applic. \textbf{41} (1941), no.~3--4,
  327--352.

\bibitem{LeanIJNM21}
Meng~H. Lean and A.~Wexler, \emph{Accurate numerical integration of singular
  boundary element kernels over boundaries with curvature}, Int. J. Num. Meth.
  Engin. \textbf{21} (1985), no.~2, 211--228. \MR{0784707}

\bibitem{Luke1975}
Yudell~L. Luke, Bing~Yuan Ting, and Marilyn~J. Kemp, \emph{On generalized
  {G}aussian quadrature}, Math.\ Comp. \textbf{29} (1975), no.~132, 1083--1093.
  \MR{0388740 (52 \#9574)}

\bibitem{SalzerCJ15}
Herbert~E. Salzer, \emph{Lagrangian interpolation at the {C}hebyshev points
  $x_{n,\nu}=\cos(\nu\pi/n)$, $\nu=0 (1) n$; some unnoted advantages}, Comp. J.
  \textbf{15} (1972), no.~2, 156--159. \MR{0315865 (47 \#4414)}

\bibitem{WynnQJM18}
P.~Wynn, \emph{A general system of orthogonal polynomials}, Quart. J. Math.
  Oxford \textbf{18} (1967), no.~1, 81--96. \MR{0210963 (35 \#1848)}

\end{thebibliography}

\end{document}